\documentclass[12pt]{article}

\begin{document}
\newtheorem{theorem}{Theorem}[section] 
\newtheorem{proposition}[theorem]{Proposition}
\newtheorem{lemma}[theorem]{Lemma}
\newtheorem{corollary}[theorem]{Corollary}
\newtheorem{Definition}[theorem]{Definition}
\newenvironment{definition}{\begin{Definition}\rm}{\end{Definition}}
\newtheorem{Remark}[theorem]{Remark}
\newenvironment{remark}{\begin{Remark}\rm}{\end{Remark}}
\newtheorem{Example}[theorem]{Example}
\newenvironment{example}{\begin{Example}\rm}{\end{Example}}
\newcommand{\Proof}{\noindent{\bf Proof:} }
\catcode`\@=11
\@addtoreset{equation}{section}
\catcode`\@=12
\renewcommand\theequation{\arabic{section}.\arabic{equation}}

\newcommand\sqr[2]{{\vcenter{\vbox{\hrule height.#2pt\hbox{\vrule width.#2pt
  height#1pt \kern#1pt \vrule width.#2pt}\hrule height.#2pt}}}}
\newcommand\square{\mathchoice\sqr64\sqr64\sqr{2.1}3\sqr{1.5}3}
\newcommand\qed{\hfill$\square$}
\newcommand\cove{{\cal A}}
\newcommand\id{{\rm id}}
\newcommand\infigamma{{\cal I}_\gamma^{\infty}}
\newcommand\II{{\cal I}}
\newcommand\FF{{\cal F}}
\newcommand\al{\alpha}
\newcommand\be{\beta}
\newcommand\ga{\gamma}
\newcommand\de{\delta}
\newcommand\om{\omega}
\newcommand\Om{\Omega}
\newcommand\iy{\infty}
\newcommand\intiy{\int_{-\ga\cdot\iy}^{\ga\cdot\iy}}
\newcommand\pa{\partial}
\newcommand\som{{\,{\rm s}\om}}
\newcommand\RR{{\bf R}}
\newcommand\CC{{\bf C}}
\newcommand\ZZ{{\bf Z}}
\newcommand\Zplus{\ZZ_{\ge0}}
\newcommand\HH{{\cal H}}
\newcommand\HD{\HH^{\rm D}}
\newcommand\HS{\HH^{\rm S}}
\renewcommand\Re{{\rm Re}}
\renewcommand\Im{{\rm Im}}
\newcommand\half{{1\over 2}}
\newcommand\mPP{\medbreak}
\newcommand\sLP{\smallbreak\noindent}
\newcommand\mLP{\medbreak\noindent}
\newcommand\bLP{\bigbreak\noindent}
\newcommand\qchoose[2]{{#1 \brack #2}_q}

\title{A $q$-analogue of convolution on the line}
\label{convolution}
\author{G. Carnovale
\thanks{Research for this paper was started while the first author worked at
the
University of Utrecht, supported by NWO, project number 610-06-100.
The paper was completed while she worked at the University of Cergy-Pontoise,
supported by the University of Trieste.}
\and \hspace{-.9 cm}and T. H. Koornwinder}

\date{July 19, 2000}

\maketitle

\begin{abstract} \noindent
In this paper we study a $q$-analogue of the convolution product on the  line
in detail. 
A convolution product on the braided line was defined
algebraically  by Kempf and Majid. We adapt their definition in order to give
an analytic definition for the $q$-convolution and we study convergence
extensively.
Since the braided line is commutative as an algebra, all results
can be viewed both as results in classical $q$-analysis and in braided
algebra.
We define various classes of functions on which the convolution is
well-defined and we show that they are algebras under the defined product.
One particularly nice family of algebras, a decreasing chain depending
on a parameter running through $(0,1]$, turns out to have $1/2$ as the
critical parameter value above which the algebras are commutative.
Morerover, the commutative algebras in this family are precisely the algebras
in which each function is determined by its $q$-moments.\\
We also treat the relationship between $q$-convolution and $q$-Fourier
transform. Finally, in the Appendix, we show an
equivalence between the existence of an analytic continuation of  a function
defined on a $q$-lattice, and the behaviour of its $q$-derivatives.   
\end{abstract}
\bLP\bLP
{\sl 1991 Mathematics Subject Classification}:
33D80, 33D15, 42A85 (primary);\\
17B37 (secondary).
\bLP
{\sl Key words and phrases}:
$q$-convolution, braided line, $q$-analysis, $q$-moments,
$q$-Gaussian,
$q$-Fourier transform, classes of holomorphic functions.
\newpage
\section{Introduction}\label{intro}
The classical Fourier transform $\FF$ and the classical convolution
product are closely tied to each other by the homomorphism property
$\FF (f*g)=\FF (f) \FF (g)$, while both operations
have a conceptual interpretation on the real line $\bf R$ considered
as locally compact abelian group.
Classical convolution is of great importance, both in theory and in
applications. Whenever one has an interesting generalization
of the classical Fourier transform (like the $q$-Fourier transform
in the present paper), it is therefore natural to consider an
analogue of classical convolution which is related to the generalized
Fourier transform in a similar way as classical convolution is
related to classical Fourier transform.
The present paper and its sequel will consider convolution related
to a $q$-Fourier transform involving the $q$-exponential
function $E_q$ as a kernel.
This transform was analytically introduced 
by the second author in \cite{Koo}, Section 8,
and it was earlier considered algebraically
by Kempf and Majid \cite{KeMa}, where it occurs as the special case for
the braided line of their general theory of Fourier transform
on ``braided covector algebras''.
The braided line, a deformation as a braided
group of the algebra of functions on $\bf R$, is the simplest
non-trivial example of a braided covector algebra.
It was first introduced in a rudimental way by the second author in Section
6.8 of \cite{Koor} 
and in full detail by S. Majid, 
see \cite{Majico}, \cite{Maj1}.
The Fourier transform on a braided covector algebra of type $A_n$
(a braided analogue of function space on ${\bf R}^n$) was
studied in more detail by the first author in \cite{Ca2}.

In \cite{KeMa} Kempf and Majid also defined
convolution for those braided co\-vector algebras which have an integral which
is
bosonic and invariant under translation.
A slight adaptation of their definition in case of the braided line
(where the integral is not bosonic) is the starting point of our
analytic definition of $q$-convolution $f*_\ga g$ on ${\bf R}$ given
below (here $\ga$ denotes the choice of a $q$-lattice).
Two other motivations for this definition can be given:
the formal limit for $q\to 1$ yields classical convolution,
and $\FF_\ga(f*_\ga g)=(\FF_\ga f)\,(\FF_\ga g)$ holds
if we take for $\FF_\ga$ the $q$-Fourier transform of \cite{Koo} involving
$E_q$ as a kernel.

The aim of the present paper is to find
suitable function spaces for $f$ and $g$
such that their $q$-convolution $f*_\ga g$ is well-defined, and to
find function classes on which associativity and commutativity hold.
Let us introduce some notation for explaining this.

Throughout this paper $q$ is fixed such that $0<q<1$.
The $q$-derivative of a function $f$ is given by
$(\pa f)(x):={{f(x)-f(qx)}\over{(1-q)x}}$, and its $q$-shift by
$(Qf)(x):=f(qx)$.
For $\ga>0$ we have the $q$-lattice
$L(\ga):=\{\pm q^{k}\ga\,|\,k\in\ZZ\}$.
For a function $f$ on $L(\ga)$ the (unbounded) $q$-integral over $L(\ga)$
is denoted and defined by
\begin{equation}\label{q-int}
\int_\ga f=
\intiy f(t)\,d_qt:=
(1-q)\sum_{k=-\iy}^\iy\sum_{\epsilon=\pm1}q^{k}\ga\,f(\epsilon q^k\ga),
\end{equation}
provided the summation absolutely converges. 

Fix $\ga>0$.
Let $\II_\ga$ denote the space of absolutely $q$-integrable
functions on $L(\ga)$, and let
$\II_\ga^\iy$ denote the subspace of functions $f\in\II_\ga$
such that  $x\mapsto f(x)x^e$ is in $\II_\ga$ for every $e\in\Zplus$.
For $f\in\infigamma$ define the {\em moments}, respectively
{\em strict moments} of $f$ by:
\begin{equation}\label{momenteq}
\mu_{e,\ga}(f):=q^{{e^2+e}\over2}\intiy f(x)x^e\,d_qx,
\quad
\nu_{e,\ga}(f):=q^{{e^2+e}\over2}\intiy |f(x)x^e|\,d_qx.
\end{equation}
The factor $q^{{e^2+e}\over2}$ is a normalization factor by means of which
most statements and formulas can be expressed in a simplified form
(see Section 4).
\begin{definition}\label{convodefi}
Let $f\in\infigamma$ and let $g$ be
a function defined on some subset of $\CC$.
Then the {\em $q$-convolution product} $f*_\ga g$ is the function
given by
\begin{equation}\label{convoeq}
(f*_{\ga}g)(x):=\sum_{e=0}^\iy{{(-1)^e
\mu_{e,\ga}(f)}\over{[e]_{q}!}}(\pa^eg)(x)
\end{equation}
for $x\in\CC$ such that the $q$-derivatives $(\pa^eg)(x)$ are well-defined
for all $e\in\Zplus$ and the
sum on the right converges absolutely,
and with notation (\ref{momenteq}) being used.\qed
\end{definition}

By the asymmetric form of our definition of $f*_\ga g$,
the initial choices of function classes for $f$ and $g$ are quite different:
for $f$ the moments $\mu_{e,\ga}(f)$
should behave as $O(q^{{\al e^2}\over2}b^e)$ for $e\to\iy$ with
$\al,b>0$ (we call this of left type $\al$),
while for $g$ the $q$-derivatives $(\pa^kg)(x)$ at $x$
should behave as $O(R^k)$ for $k\to\iy$ with $R>0$.
We also need functions on $L(\ga)$ which extend to analytic functions,
on a disk centered at 0, or on a strip around $\bf R$.
Some equivalent characterizations of function spaces derived in this paper
may have independent interest (see for instance the Appendix).

Commutativity of the convolution product
is the hardest and most interesting issue of this paper.
Both the homomorphism property and the interpretation on the
(commutative) braided line suggest commutativity, but we find significant
counterexamples. An explanation is that the $q$-Fourier transform of $f$
only depends on the moments of $\mu_{e,\ga}(f)$ and that
$\mu_{e,\ga}(f*_\ga g)$ is symmetric in $f$ and $g$, but that
$f$ in certain function classes is not completely determined by its moments.
For getting commutativity we need functions $f$ (and $g$) in the
convolution product for which the strict moments
$\nu_{e,\ga}(f)$
behave as $O(q^{{\al e^2}\over2}b^e)$ for $e\to\iy$ with
$\al>1/2$, $b>0$ (we call this of strict left type $\al>1/2$),
and which are holomorphic on
a strip around $\bf R$.
It turns out that the $q$-Gaussian $x\mapsto e_{q^2}(-x^2)$,
which has has strict left
type $1/2$,
does not commute with many entire functions of strict left type $>1/2$.

Part of the results of this paper occurred in the recent dissertation
\cite{giova} by the first author.
She discusses some further aspects of $q$-convolution
in the subsequent paper \cite{Ca3}.
\mLP
{\bf Acknowledgement}\qquad
We thank the referee for his careful reading of the manuscript and
his constructive remarks.
\mLP
{\bf Further notations}\qquad
We denote as usual:
\begin{equation} \label{usual}
{(a;q)_k:=\prod_{j=0}^{k-1}(1-aq^{j}),\quad
(a;q)_{\iy}:=\lim_{k\to\iy}(a;q)_k\,,
\atop
[k]_{q}:={{1-q^{k}}\over{1-q}}\,,\quad
[k]_{q}!:={{(q;q)_k}\over{(1-q)^k}}\,,\quad
\qchoose k j:=
{{[k]_q!}\over{[j]_q![k-j]_q!}}={{(q;q)_k}\over{(q;q)_j(q;q)_{k-j}}}\,.}
\end{equation}
\mLP
A function $f\colon x\mapsto f(x)$ may also be denoted as $f(X)$.
This will be useful for functions like
$fX^e\colon x\mapsto f(x)x^e$ and
$e_{q^2}(-X^2)\colon x\mapsto e_{q^2}(-x^2)$.\\
For $q$-hypergeometric series the notation of
Gasper \& Rahman \cite{gasper} will be followed.\\
The (bounded) $q$-integral of a function $f$ on 
$\{\pm\ga q^k\mid k\in\Zplus\}$
is denoted and defined by
\begin{equation}
\int_{-\ga}^\ga f(t)d_qt=(1-q)\sum_{k=0}^\iy
\sum_{\epsilon=\pm1} q^k\ga f(\epsilon q^k \ga)
\end{equation}
provided the summation absolutely converges.
\section{Motivation of the definition of $q$-convolution}
We will give three different motivations for our Definition
\ref{convodefi} of $q$-convolution
in the following three remarks.
\begin{remark}\label{braidremark}
The {\em braided line} (see \cite{Maj1})
is a braided Hopf algebra $\cove$ which, as an algebra, is equal to the
(commutative) algebra $\CC[[x]]$ of formal power series
in $x$, and which has
braiding $\Phi(x^k\otimes x^l):=q^{kl} x^l\otimes x^k$,
comultiplication
$\Delta(x^k):=\sum_{j=0}^k\;\qchoose k j x^{k-j}\otimes x^j$,
counit
$\varepsilon(x^k):=\de_{k,0}$
and braided antipode
$S(x^k):=(-1)^kq^{{k^2-k}\over2}x^k$.
Then the $q$-analogue of Taylor's formula is given by
\begin{equation}\label{q-Taylor}
\Delta(f(x))=\sum_{j=0}^\iy{{x^j}\over{[j]_{q}!}}\otimes
\pa^j(f(x)),
\end{equation}
where $\pa$ denotes the $q$-derivative.

The original formal definition of a convolution on any braided covector
algebra $A$ with a bosonic integral 
$\int\colon A\to \CC$ invariant under translation, was given by Kempf \&
Majid \cite{KeMa} as follows:
\begin{equation}\label{Majid-convo}
(f*g)(x):=({\textstyle\int}\otimes\id)(m\otimes\id)(\id\otimes
S\otimes\id)(\id\otimes \Delta)(f\otimes g)(x).
\end{equation}
In the case of the braided line we take for $\int$
the integral defined by (\ref{q-int}). This integral is
invariant under translation (see \cite{KeMa}, and \cite{Koo} and \cite{Ca2}
for an analytic proof), but it is not bosonic.
Associativity will fail for the convolution defined by (\ref{Majid-convo}).
Therefore, we slightly modify (\ref{Majid-convo}) into
\begin{equation}\label{our-convo}
(f*_{\ga}g)(x):=
({\textstyle\int}\otimes\id)
(m\otimes\id)
(\id\otimes Q\otimes\id)
(\id\otimes S\otimes\id)
(\id\otimes \Delta)
(f\otimes g)(x).
\end{equation}
When we substitute the $q$-Taylor formula (\ref{q-Taylor})
into (\ref{our-convo}) then we
formally get formula (\ref{convoeq}) (with (\ref{momenteq}) substituted),
i.e., our original definition of $q$-convolution.

We  may interpret $f*_{\ga}g$ formally as
the action of a pseudo-$q$-differential operator ${\cal D}f(\pa)$ on
$g(x)$, where ${\cal D} f(\pa):=\FF '_S(\ga)(f)(i(1-q)\pa)$.
Here
$\FF '_S(\ga)(f):=\sum_{e=0}^\iy{{(-1)^e\mu_{e,\ga}(f)}\over{[e]_{q}!}}
\otimes\pa^e$  is the {\em weak braided Fourier transform} considered in
\cite{Ca2} (with $q^2$ instead of $q$).
\hfill$\spadesuit$\end{remark}
\begin{remark}
The formal limit for $q\uparrow 1$ of formula (\ref{convoeq}) is
the classical convolution product:
\begin{eqnarray*}
\lefteqn{\sum_{e=0}^\iy\biggl(\int_{-\iy}^\iy dt\, f(t){{(-1)^et^e}\over{e!}}
\biggr)g^{(e)}(x)\mathrel{\mathop=^{\rm formally}}
\int_{-\iy}^\iy dt \Biggl(f(t)
\biggl(\sum_{e=0}^\iy{{(-1)^et^e}\over{e!}}g^{(e)}(x)\biggr)\Biggr)}
\hspace{2in}\\
&&\mathrel{\mathop=^{\rm formally}} \int_{-\iy}^\iy dt\, f(t)g(x-t)=(f*g)(x).
\hspace{.53in}\spadesuit
\end{eqnarray*}
\end{remark}
\begin{remark}\label{fourier-remark}
Recall the two $q$-exponentials (see \cite{gasper}):
\begin{equation}\label{q-exponential}
e_{q}(x):=\sum_{k=0}^\iy{{x^k}\over{(q;q)_k}}={1\over (x;q)_\iy},\quad
E_q(x):=\sum_{k=0}^\iy{{q^{k(k-1)\over 2}x^k}\over{(q;q)_k}}=(-x;q)_\iy,
\end{equation}
where $|x|<1$ in the infinite sum defining $e_q(x)$.
In \cite{Koo} a $q$-Fourier transform pair was presented as folllows:
\begin{equation}\label{q-fourier}
\psi(x)={1\over c_q(\ga) b_q}\int_{-1}^1 e_q(-ixy)\,\phi(y)\,d_qy,
\quad
\phi(y)=
\int_{-\ga.\iy}^{\ga.\iy}E_q(iqxy)\,\psi(x)\,d_qx,
\end{equation}
where $b_q$ and $c_q(\ga)$ are given by
\begin{eqnarray} 
&&\;\quad b_q:=\int_{-1}^1 E_{q^2}(-q^2x^2)\,d_qx=
(1-q)\,(q,-q,-1;q)_\iy,\label{biqu}\\
&&c_q(\ga):=\int_{-\ga.\iy}^{\ga.\iy} e_{q^2}(-x^2)\,d_qx
={2(1-q)\,(q^2,-q\ga^2,-q\ga^{-2};q^2)_\iy\,\ga\over
(-\ga^2,-q^2/\ga^2,q;q^2)_\iy}\,.\label{ciqu}
\end{eqnarray}
Write the second transform in (\ref{q-fourier}) as $\phi=\FF_\ga\,\psi$.
Then an immediate formal computation shows that
\begin{equation}\label{homomor}
(\FF_\ga(f*_\ga g))(x)=
(\FF_\ga f)(x)\, (\FF _\ga g)(x).
\end{equation}
The transform $\FF_\ga$ is essentially the weak braided Fourier
transform $\FF'_S(\ga)$ (see Remark \ref{braidremark} and
reference \cite{Ca2}).
Equation (\ref{homomor}) will be rigorously proved for suitable $f$ and $g$ in
Section \ref{convofourier}.
\hfill$\spadesuit$
\end{remark}

For later use we recall the formulas (9.8), (9.14) in \cite{Koo}
(for $n\in\Zplus$):
\begin{eqnarray}
\label{eq-moment}
\intiy x^n\,e_{q^2}(-x^2)\,d_qx=c_q(\ga)\,q^{-k^2}\,(q;q^2)_k
\quad\hbox{if $n=2k$\quad($=0$ otherwise)}\\
\label{Eq-moment}
\int_{-1\cdot\iy}^{1\cdot\iy}
x^n\,E_{q^2}(-x^2)\,d_qx=b_q\,q^{2k+1}\,(q;q^2)_k
\quad\hbox{if $n=2k$\quad($=0$ otherwise).}\quad
\end{eqnarray}
\section{Good function spaces for $q$-convolution}\label{convoconve}
In general, it is not true that $(f*_\ga g)(x)$ can be expanded as a
(possibly formal) power series if $g(x)$ has a convergent power series
expansion. The reason is that the coefficients of an
expansion in powers of $x$ of
$(f*_\ga g)(x)$ will be in general infinite series themselves. We want to find
conditions on $f$ and $g$ so that, with the given definition,
$(f*_{\ga}g)(x)$ makes sense on some subset of the complex plane. One should
give conditions on the growth of the moments $|\mu_{e,\ga}(f)|$ of $f$
and of the power series coefficients of $g$.
For this purpose we will now introduce the class of functions
of left type $\al$.
For convenience, we also give here some variants of this definition
and corresponding notation which will be needed later in the paper.
Recall the definitions of the function spaces $\II_\ga$ and $\II_\ga^\iy$,
given in the Introduction.

\begin{definition} \label{left} \quad
\sLP
{\bf (a)}
For $\al>0$ the space $\II_{\ga,\al}^\om$ of functions
{\em of left type} $\al$ on $L(\ga)$ consists of
all $f\in \infigamma$
such that, for some $b>0$,
$|\mu_{e,\ga}(f)|=O(q^{{\al e^2}\over2}b^e)$
as $e\to\iy$.\\
The space $\II_\ga^\om$ consists of all functions of some left type $\al>0$
on $L(\ga)$.
\sLP
{\bf (b)}
For $\al>0$ the space $\II_{\ga,\al}^\som$ of functions
{\em of strict left type} $\al$ on $L(\ga)$ consists of
all $f\in\infigamma$
such that, for some $b>0$,
$\nu_{e,\ga}(f)=O(q^{{\al e^2}\over2}b^e)$
as $e\to\iy$.\\
The space $\II_\ga^\som$ consists of all functions of some strict
left type $\al>0$
on $L(\ga)$.
\sLP
{\bf (c)}
For $a>0$ the space $\HD_a$ consists of all functions which are holomorphic
on the disk $\{z\in\CC\mid |z|<a\}$.
By $\HD$ we denote the space of all functions which are holomorphic
on some disk centered at 0.
\sLP
{\bf (d)}
For $a>0$ the space $\HS_a$ consists of all functions which are holomorphic
on the strip $\{z\in\CC\mid |\Im(z)|<a\}$.
By $\HS$ we denote the space of all functions which are holomorphic
on some strip around $\bf R$.\qed 
\end{definition}
\sLP
The ``intersection'' of a $\HH$-space defined in (c) or (d) with one
of the $\II$-spaces (i.e. $\II_\ga$,
$\infigamma$ or a space defined
in (a) or (b)) will be denoted by putting the two
symbols $\HH$ and $\II$ behind each other.
Here we mean intersection in a special sense.
For instance, $\HD\infigamma$ will denote the space of the functions on
$L(\ga)$ belonging to $\infigamma$ which coincide within some disk centered
at 0 with the restriction of a (necessarily unique) holomorphic
function on that disk.
Note that, conversely, a function $f\in \HD\infigamma$ is in general
not uniquely determined
by the function in $\HD$ with which it has a common restriction
within some disk.\\
We will always assume for the ``intersections'' just defined that $\ga$ is
less than the parameter $a$ (radius of a disk or half width of
a strip) occurring in
$\HD_a$ or $\HS_a$. Clearly, this assumption is not restrictive,
since $\ga$ may be replaced by $q^k\ga$ for arbitrary $k\in\Zplus$.

Observe that, if $f\in\infigamma$ is of left type $\al$ 
(respectively strict left type $\al$) on $L(\ga)$ then
so is $X^k f$ for for each $k\in\Zplus$. Indeed,
$\mu_{e,\ga}(fX^k)=q^{-{k^2\over 2}-{k\over 2}-ek}\mu_{k+e,\ga}(f)$,
and similarly for $\mu$ replaced by $\nu$.
\begin{example}\label{Gaussian} \quad
\sLP
{\bf (a)}
By (\ref{eq-moment}) the $q$-Gaussian $e_{q^2}(-X^2)$ belongs to $\infigamma$
for each $\ga>0$ and
$$
\mu_{2k,\ga}(e_{q^2}(-X^2))=c_q(\ga)\,(q;q^2)_k\,q^{k^2+k},\quad
\mu_{2k+1,\ga}(e_{q^2}(-X^2))=0.
$$
Hence $e_{q^2}(-X^2)$ is of left type $1/2$ on $L(\ga)$.
It is also of strict left type $1/2$ on $L(\ga)$.
Indeed, if $f\in\infigamma$
is an even and nonnegative function then
$\nu_{2k,\ga}(f)=\mu_{2k,\ga}(f)$ and
$\nu_{2k+1,\ga}(f)\le \mu_{2k+2,\ga}(f)+\mu_{2k,\ga}(f)$.
\sLP
{\bf (b)}
By (\ref{Eq-moment}) the $q$-Gaussian $E_{q^2}(-X^2)$
belongs to $\II_1^\iy$ and
$$
\mu_{2k,1}(E_{q^2}(-X^2))=b_q\,(q;q^2)_k\,q^{2k^2+k},\quad
\mu_{2k+1,1}(E_{q^2}(-X^2))=0.
$$
Hence $E_{q^2}(-X^2)$ is of left type $1$ on $L(1)$.
It is also of strict left type $1$ on $L(1)$ by a similar argument as in (a).
\hfill$\spadesuit$
\end{example}
\begin{remark}
Let $g\in\HD_R$ with $g(x)=\sum_lc_lx^l$  and put
$G(r):=\sum_l|c_l|r^l$
(we will repeatedly use this convention).
Note that then also $\pa^k g\in \HD_R$ for every $k\in\Zplus$.
Since $|c_l|\le\bigl(\sum_p|c_p|r^p\bigr)r^{-l}$ for every $r<R$ and since
$(q^l;q^{-1})_k\le1$
for $l\ge k$ we can also say that, for every $|x|<r<R$,
\begin{equation}\label{gG}
|(\pa^kg)(x)|\le
{{r}\over{r-|x|}}\,{1\over{r^k(1-q)^k}}\,G(r).
\end{equation}
It will be shown in the Appendix that condition (\ref{gG}) together with
equality
of left and right $q$-derivatives at $0$ is equivalent to analyticity.
\hfill$\spadesuit$
\end{remark}
\begin{lemma}\label{welldefined}
Let $f\in\II_{\ga,\al}^\om$ and $g\in \HD_a$. 
Then $f*_{\ga}g\in \HD_a$.
If $g$ can be continued analytically on a starlike
domain $\Om$ centered at $0$, then $f*_{\ga}g$ is well-defined and
analytic on $\Om$.
\end{lemma}
\Proof
We want to show that $(f*_{\ga}g)(x)$
has a convergent power series expansion for $|x|<a$.
Let $g(x):=\sum_lc_lx^l$.
Then direct computation gives
\begin{equation}\label{doublesum}
(f*_{\ga}g)(x)=\sum_{e=0}^{\iy}\sum_{l=e}^{\iy}
{{(-1)^e\mu_{e,\ga}(f)}\over{[e]_{q}!}}\,c_l\,{{[l]_{q}!}
\over{[l-e]_{q}!}}\,x^{l-e}.
\end{equation}
For every $r<a$ and every $x\in\CC$ such that $|x|<r$ one has
\begin{eqnarray*}
\lefteqn{\sum_{e=0}^{\iy}\sum_{l=e}^{\iy}\left|\mu_{e,\ga}(f)\right|\,|c_l|
\qchoose l e  |x|^{l-e}
\le  
C\sum_{e=0}^{\iy}\sum_{p=0}^{\iy}{{q^{{{\al}\over2}
e^2}b^e(q;q)_{e+p}}\over{(q;q)_p(q;q)_e}}\,G(r)r^{-p-e}|x|^p}\hspace{1.75in}\\
&&\le G(r)\,C \sum_{e=0}^{\iy}
{{q^{{{\al}\over2}e^2}r^{-e}b^e}\over{(q;q)_e}}
\sum_{p=0}^{\iy}\left({{|x|}\over{r}}\right)^p <\iy
\end{eqnarray*}
for some $C>0$. 
Hence, by dominated convergence, we may invert the order of summation in
formula (\ref{doublesum}) so that, for $|x|<a$, we have
$$
(f*_{\ga}g)(x)=\sum_{p=0}^{\iy}
\left(\sum_{e=0}^{\iy}(-1)^e\mu_{e,\ga}(f)
\qchoose{p+e} e c_{p+e}\right)x^p=
\sum_{p=0}^{\iy}b_px^p
$$
where $|b_p|\le G(r)\,C\,r^{-p}\sum_{e=0}^\iy{{q^{{\al\over2}
e^2}(q^{p+1};q)_e}\over{(q;q)_e}}\,b^er^{-e}<\iy$, while the power
series of $(f*_{\ga}g)(x)$ converges absolutely for $|x|<a$.

Let now $g$ be continued analytically on a starlike domain $\Om$
centered at $0$. 
It follows by induction that for every starlike compact
set $K$ centered at $0$ and contained
in $\Om$, $|(\pa^eg)(x))|\le ||g||_K\,{{2^e}\over{|x|^e(1-q)^e}}$ for every
nonzero $x\in K$, where $\parallel\cdot\parallel_{K}$ denotes the supremum
norm. Let $s>0$. Then, for $x\in K$ with $|x|\ge s$,
$$
(f*_{\ga}g)(x)=\sum_{e=0}^\iy{{(-1)^e\mu_{e,\ga}(f)}\over{[e]_{q}!}}
(\pa^eg)(x)
$$
with  
$$
\biggl|{{(-1)^e\mu_{e,\ga}(f)}\over{[e]_{q}!}}(\pa^eg)(x)\biggr|\le
C{{q^{{\al\over2} e^2}b^e}\over{[e]_{q}!}}|(\pa^eg)(x)|\le
C{{q^{{\al\over2} e^2}(2b)^e}\over{s^e(q;q)_e}}\,.
$$
Hence, on each set
$\{x\in K\,|\,|x|>s\}$ with $K\subset \Om$ a starlike compact set centered at
$0$ and $s>0$, the series expressing
$(f*_{\ga}g)(x)$ is uniformly convergent, with\\
$\sum_{e=0}^E{{(-1)^e\mu_{e,\ga}(f)}\over{[e]_{q}!}}(\pa^eg)(x)$ analytic
in $x$ for every
$E\in\Zplus$. It follows that $f*_{\ga}g$ is analytic on each set $\{x\in
\Om\mid|x|>s\}$ with $s>0$. Since we already showed that $f*_{\ga}g$ is
analytic on $\{x\in\CC\mid |x|<a\}$, this completes the proof.
\qed
\mPP
Now we consider $g$ in $\HD_a\II_\ga$ or
in $\HS_a\II_\ga$ with $a>\ga$.
An example of a
function  in $\HD_a\II_\ga$ is any function $g$ with finite support
contained in $\{\pm q^{-k}\ga\,|\,k\in\ZZ_{\ge1}\}$. 

An example of a function in $\HS_a\II_\ga$ is
the function $e_{q^2}(-X^2)$. It is clearly an element of  $\HS_1\II_\ga$
for every $\ga\in(0,\,1)$, since
$e_{q^2}(-x^2)=\sum_{k=0}^\iy{{(-1)^kx^{2k}}\over{(q^2;q^2)_k}}$
for $|x|<1$ and since this function has analytic continuation
$1\over{(-x^2;q^2)_\iy}$ for $x\ne iq^{-k}$ ($k\in\Zplus$).
\begin{lemma}\label{integrapartial}
Let $g\in\HD_a\II_\ga^\iy$ for some $a>\ga$.
Then for every $r\in(\ga,\,a)$  there is a constant $B>0$ such that 
$$
\intiy|x^k(\pa^eg)(x)|\,d_qx \le
\Bigl(\intiy|x^kg(x)|\,d_qx+r^k{B}\Bigr){{2^e}\over{\ga^e(1-q)^e}}
$$
for every $e$ and $k$ in $\Zplus$. 
Hence $\pa^e g\in\HD_a\II_\ga^\iy$.\\
\noindent If $g\in \HD_a\II_\ga$ but not necessarily in $\II_\ga^\iy$, 
the above conclusion still holds for $k=0$ and $e\in\Zplus$. Then
$\pa^e g\in\HD_a\II_\ga$.
\end{lemma} 
\Proof
If $e=0$ the statement is trivial. Let $e>0$ and
$g\in\HD_a\II_\ga^\iy$. Then
for every $r\in(\ga,\,a)$ one has:
\begin{eqnarray*}
&&\intiy|x|^k|(\pa^eg)(x)|\,d_qx=
\int_{-\ga}^{\ga}|x^k(\pa^eg)(x)|d_{q}x\\ 
&&+(1-q)\sum_{s\le
-1}\sum_{\epsilon=\pm1}q^{s(k+1)}\ga^{k+1} {{|(\pa^{e-1}g)(\epsilon
q^{s}\ga)-(\pa^{e-1}g)(\epsilon
q^{s+1}\ga)|}\over{(1-q)q^{s}\ga}}\\ 
&&\le \sum_{s\le-1}\sum_{\epsilon=\pm1}q^{ks}
\ga^k\biggl[\bigl|(\pa^{e-1}g)(\epsilon
q^{s}\ga)\bigr|+\bigl|(\pa^{e-1}g)(\epsilon
q^{s+1}\ga)\bigr|\biggr]+{{2r^{k-e}G(r)\,r\ga}\over{(r-\ga)(1-q)^e}}\\ 
&&\le{2\over{\ga(1-q)}}\intiy|x^k(\pa^{e-1}g)(x)|\,d_qx+{{2r^{k}
\ga^{-e}
G(r)\,r\ga}\over{(r-\ga)(1-q)^{e}}}
\end{eqnarray*} where the first inequality follows by estimate (\ref{gG}).
Hence, by iterating the process we get:
\begin{eqnarray*}
&&\intiy|x^k(\pa^eg)(x)|\,d_qx \le {{2^2}\over{\ga^2(1-q)^2}}
\intiy |x^k(\pa^{e-2}g)(x)|\,d_qx\\
&&\hspace{1.75in}+{{2r^{k}\ga^{-e}
G(r)\,r\ga}\over{(r-\ga)(1-q)^{e}}}
+{{2^2r^{k}\ga^{-e+1-1}
G(r)\,r\ga}\over{(r-\ga)(1-q)^{e-1+1}}}\le\ldots\\
&&\le {{2^e}\over{\ga^e(1-q)^e}}\intiy |x^k(g(x))|\,d_qx+{{2^er^{k}
\ga^{-e}
G(r)\,r\ga}\over{(r-\ga)(1-q)^{e}}}\biggl[1+{1\over2}+
\cdots\biggl({1\over2}\biggr)^{e-1}\biggr].
\end{eqnarray*}
Therefore we get the statement with $B={{2
G(r)r\ga}\over{(r-\ga)}}$. 
The conclusion for $k=0$ if $g\in\HD_a\II_\ga$ is also clear from the proof.
\qed
\mLP
\begin{remark}
If $g\in\HS_a\II_\ga$, then
$\int\pa^kg(x)\sim_{\ga}0$ and the $q$-integral of $g$ is {\em always}
invariant under translation, in the notation of Section IV in \cite{Ca2}. In
particular, this result will then follow without assuming
condition (c) in that Section, since by the above Lemma the
partial $q$-derivatives of $g(x)$ are
automatically $q$-integrable.\hfill$\spadesuit$
\end{remark}
\smallbreak
{}From now on we will often use the shorthand notation $\int_\ga f$ for
$\intiy f(t)\,d_qt$, see formula (\ref{q-int}).
We can conclude:
\begin{proposition}\label{fconvogintegrable}
Let $f\in\II_{\ga'}^\om$ and 
$g\in\HD_a\II_\ga$. Then $f*_{\ga'}g$ is well-defined
and absolutely $q$-integrable on $L(\ga)$, and $f*_{\ga'}g\in\HD_a\II_\ga$.
\end{proposition}
\Proof
Let $d:=2/(\ga(1-q))$.
By Lemma \ref{integrapartial} we know that 
$\int_\ga|\pa^eg|\le Ad^e$ for some
$A>0$. 
Hence for some $B>0$ we have that
$|(\pa^eg)(\pm q^{-k}\ga)|\le B d^e q^{k}$
for every
$e,\,k\in\Zplus\;$.
Then the proof that $f*_{\ga'}g$ is well-defined on $L(\ga)$
is similar to the second part of the proof of
Lemma \ref{welldefined}.

We still have to show that $\int_\ga|f*_{\ga'}g|<\iy$. 
By assumption,
$|\mu_{e,\ga'}(f)|\le Cq^{{\al\over2} e^2}b^e$ for some constants
$C,\,b$ and
$\al>0$. Then
\begin{eqnarray*}
&&\int_\ga \,\biggl|\sum_{e=0}^\iy{{(-1)^e\mu_{e,\ga'}(f)}
\over{[e]_{q}!}}\pa^eg\biggr|\\
&&\qquad\le
\sum_{e=0}^\iy{{|\mu_{e,\ga'}(f)|}\over{[e]_{q}!}}
\int_\ga \left|\pa^eg\right|
\le
{{AC}\over{(q;q)_\iy}}
\sum_{e=0}^\iy q^{{{\al}\over2}
e^2}{d^eb^e(1-q)^e}<\iy.
\quad\square
\end{eqnarray*}
\begin{corollary} Let $f\in\II_{\ga'}^\om$ and 
let $g\in\HD_a\II_\ga$. Then
$\pa^k(f*_{\ga'}g)$ and $f*_{\ga'}\pa^kg$
are also absolutely $q$-integrable on
$L(\ga)$ for every $k\in\Zplus$.\end{corollary}
\Proof
By Proposition \ref{fconvogintegrable} we know that 
$f*_{\ga'}g\in\HD_a\II_\ga$. Hence
$\pa^k(f*_{\ga'}g)\in\HD_a\II_\ga$ by Lemma 
\ref{integrapartial}. By Lemma \ref{integrapartial} we also know that
$\pa^kg\in\HD_a\II_\ga$, so that
$f*_{\ga'}\pa^kg\in\HD_a\II_\ga$ by Proposition \ref{fconvogintegrable}.\qed
\section{Associativity of $q$-convolution}\label{asso}
The next step will be to investigate associativity. In order to do this, we
need to know under which hypotheses  $(f*_\ga g)*_{\ga}h$ and
$f*_\ga(g*_\ga h)$ are well-defined.  For $f*_\ga(g*_\ga h)$
this was essentially described in Lemma
\ref{welldefined} and Proposition \ref{fconvogintegrable}, i.e., $f$ and
$g$ need to be of left type,
and $h$ has to belong to $\HD$.\\ 
In order to
understand when 
$(f*_\ga g)*_{\ga}h$ is well-defined, and to prove associativity,
we need to
investigate the behaviour of
$\mu_{e,\ga'}(f*_{\ga}g)$. We will use  the following lemmas.
\begin{lemma}\label{movepartial}
Let $f\in\infigamma$ such that
$\pa^j f\in\infigamma$
for every $j\in\Zplus$. Then
$$
\intiy
(\pa^af)(x)\,x^b\,d_qx=\cases{(-1)^{a}q^{-ba+{{a^2-a}\over2}}{{[b]_{q}!}
\over{[b-a]_{q}!}}
\intiy f(x)\,x^{b-a}\,d_qx& if $b\ge a$,\cr 0& otherwise.\cr}
$$
In particular,
\begin{equation}\label{momentpartial}
\mu_{e+k,\ga}(\pa^k f)=
(-1)^k{{[e+k]_q!}\over{[e]_q!}}\mu_{e,\ga}(f)
\quad\hbox{for  every $k\in\Zplus\,$,}
\end{equation}
and $\mu_{l,\ga}(\pa^k f)=0$ if $l<k$.
\end{lemma}
\Proof
It is easy to see that
$q^{b}X^b\,\pa f=\pa (f\,X^b)-f\,\pa(X^b)$ (this is a particular case of the
braided Leibniz rule in \cite{KeMa} applied to $X^b$ and $f$). If we apply
 the unbounded $q$-integral
to both sides with $f$ replaced by $\pa^{a-1}f$, then by absolute
$q$-integrability we obtain:
$$
\int_\ga X^b\,\pa^a f=-q^{-b}[b]_{q}\int_\ga 
X^{b-1}\,\pa^{a-1}f.
$$
By repeating this procedure we get the statements of the Lemma. \qed
\begin{lemma}\label{convoleftisleft}
If $f\in\HD\II_{\ga,\al}^\om$ then
$X^k\,\pa^j f\in\HD\II_{\ga,\al}^\om$ for every $k\in\Zplus$.
\end{lemma}
\Proof
By  Lemma \ref{integrapartial}
$X^k\,\pa^j f$ is absolutely $q$-integrable on
$L(\ga)$ for every $j$ and $k$.
By Lemma \ref{movepartial} $\mu_{k+e,\ga}(\pa^j f)=0$
unless $e+k>j$. For $e+k>j$ we have
\begin{eqnarray*}
\lefteqn{|\mu_{e,\ga}(X^k\,\pa^jf)|=|\mu_{k+e-j}(f)|q^{-{{k^2+2ek+k}\over2}}
{{(q^{k+e-j+1};q)_{j}}\over{(1-q)^j}}}\hspace{1in}\\
&& \le \biggl[{{Cq^{{-k^2-k}\over2}b^{k-j}}\over{(1-q)^j}}
q^{{\al\over2}(k-j)^2}\biggr]
q^{-e(k-\al k+\al j)}q^{{\al\over2} e^2}b^e
\end{eqnarray*}
for some $b>0$. 
Hence $|\mu_{e,\ga}(X^k\,\pa^jf)|=
O\bigl(q^{{\al\over2} e^2}(bq^{-(k-\al k+\al j)})^e\bigr)$
as $e\to\iy$.\qed
\mPP
It follows that $\pa^k f*_{\ga}g$ is well-defined if
$f\in\HD_a\II_\ga^\om$
and $g\in\HD$.
\begin{lemma}\label{interchange} Let
$f\in\HD_a\II_\ga^\om$ and let $g$ be defined,
together with its $q$-derivatives, 
on a domain $\Om$. Let $0\not=x\in\Om$ be such that, for some $R_x>0$,
$|\pa^eg(x)|=O(R_x^e)$ as 
$e\to\iy$. Then,
for every $k\in\Zplus$, $\pa^k(f*_{\ga}g)$, 
$(\pa^k f*_{\ga}g)$ and 
$(f*_{\ga}\pa^k g)$ are well defined at $x$ and
$$
\pa^k(f*_{\ga}g)(x)=
(\pa^k f*_{\ga}g)(x)=
(f*_{\ga}\pa^k g)(x).
$$
In particular, the result holds
if $g\in\HD\II_{\ga}$ and $x\in L(\ga)$, or
if $g$ is analytic on a starlike domain $\Om$ centered at $0$
and $x\in\Om$.\end{lemma}  
\Proof Under the hypothesis on $f$, $g$ and $x$, the convolution product
$f*_\ga g$ is well defined at $x$ by the proof of Lemma \ref{welldefined} and
Proposition \ref{fconvogintegrable}.
If $g$ and $x$ satisfy the hypothesis
then also $\pa^k g$ and $x$ satisfy the hypothesis
for every $k\in\Zplus$, hence
$f*_\ga \pa^k g$  and $\pa^k(f*_\ga g)$ are well-defined at $x$. By
Lemma
\ref{convoleftisleft} $(\pa^kf)*_\ga g$ is also well defined at $x$. The three
expressions are: 
\begin{eqnarray*}
&&(\pa^k(f*_{\ga}g))(x)=
\pa^k\left(\left(
\sum_{e=0}^\iy{{(-1)^e\mu_{e,\ga}(f)}\over{[e]_{q}!}}\right)
\pa^eg\right)(x)=(f*_{\ga}\pa^k g)(x)\;\;\,{\hbox{and}}\\ 
&&(\pa^k f*_{\ga}g)(x)=\sum_{e\ge0}
{{(-1)^{e+k}\mu_{e+k,\ga}(\pa^k f)}\over{[e+k]_{q}!}}
(\pa^{e+k}g)(x)
\qquad\qquad\qquad\qquad\qquad\\
&&\hspace{1.9in}=\sum_{e\ge0}{{(-1)^e\mu_{e,\ga}(f)}\over{[e]_{q}!}}
(\pa^{k+e}g)(x)=(f*_{\ga}\pa^kg)(x)
\end{eqnarray*}
where absolute convergence follows by the conditions imposed on $x$ and $g$,
hence the first statements.

If $g\in\HD\II_\ga$ and $x\in L(\ga)$ then
it follows by the proof of Proposition
\ref{fconvogintegrable} that $g$ and $x$ satisfy the conditions of
the Lemma.\\
If $g$ is analytic on a starlike domain $\Om$ centered at $0$ then the
condition on $g$ and
$x$ holds by the estimate in the proof of Lemma \ref{welldefined}.
Observe that in this case analytic continuation allows also $x=0\in\Om$.\qed 
\begin{proposition}\label{fede}
Let $f\in \HD\II_{\ga'}^\om$ and $g\in\HD_a\II_\ga^\iy$.
Then
$X^k\,\pa^l(f*_{\ga'}g)\in\HD_a\II_\ga^\iy$ for all $k,l\in\Zplus$.
\end{proposition}
\Proof By Proposition \ref{fconvogintegrable} $f*_{\ga'}g\in\HD_a$.
Since $\pa^l(f*_{\ga'}g)=f*_{\ga'}\pa^l g$ by Lemma
\ref{interchange} and $\pa^l g\in\HD_a\II_\ga^\iy$
by Lemma \ref{integrapartial},
we might
as well reduce to the case $l=0$. Then
\begin{eqnarray*}
\lefteqn{\int_\ga |X^k|\,|f*_{\ga'}g|=
\int_\ga |X^k|\,\biggl|\sum_{e=0}^{\iy}{{(-1)^e\mu_{e,\ga'}(f)}
\over{[e]_{q}!}}\pa^eg\biggr|}\\
&&\le C\int_\ga \sum_{e=0}^{\iy}|X^k|{{q^{{\al\over2}
e^2}b^e}\over{[e]_{q}!}}\bigl|\pa^eg\bigr|=
C\sum_{e=0}^{\iy}{{q^{{\al\over2} e^2}b^e}\over{[e]_{q}!}}
\int_\ga |X^k|\,|\pa^e g|
\end{eqnarray*}
for some $C,\,\al,b>0$.
Again by Lemma \ref{integrapartial},
there exists for any $r\in(\ga,\,a)$
a constant $B>0$ such that
$$
\int_\ga |X^k|\left|f*_{\ga'}g\right|\leq
C\sum_{e=0}^{\iy}{{q^{{\al \over2}e^2} (2b)^e}\over{(q;q)_e
\ga^e}}\Biggl(\int_{\ga'}|X^k\,f|+
r^k{B}\Biggr)<\iy.\qquad\qquad\qquad\square
$$

By the above results we know that for suitable $f$ and $g$
their convolution product 
$f*_{\ga'}g\in\HD_a\infigamma$. Next question is
then whether $f*_{\ga'}g$ is of left type too. The answer is positive.
We need the following Lemma.
\begin{lemma}\label{momenta}
Let $f\in \HD\II_\ga^\om$,
$g\in \HD\II_{\ga'}^\om$,
$k\in\Zplus$. Then
\begin{equation}\label{momentconvo}
\mu_{k,\ga'}(f*_\ga g)= 
\sum_{e=0}^k\,\qchoose k e\mu_{e,\ga}(f)\,\mu_{k-e,\ga'}(g).
\end{equation}
Thus, if $\ga=\ga'$ then
$\int_\ga (f*_{\ga}g)X^k=\int_\ga (g*_{\ga}f)X^k$ for every
$k\in\Zplus$ and $\int_\ga  (f*_\ga g)=(\int_\ga f)\,(\int_\ga g)$.
\end{lemma}
\Proof
By Proposition \ref{fede} $(f*_{\ga}g)\,X^k$ is
absolutely $q$-integrable on $L(\ga')$. Then
$$
\mu_{k,\ga'}(f*_\ga g)=q^{{k^2+k}\over2}
\int_{\ga'}(f*_{\ga}g)\,X^k=\int_{\ga'}\biggl(\sum_{e=0}^\iy
{{(-1)^e\mu_{e,\ga}(f)}\over{[e]_{q}!}}\,X^k\,\pa^e g\biggr)
$$  
Again by Lemma \ref{integrapartial} and by dominated convergence we may
interchange
integration and summation over $e$. Then, by Lemma \ref{movepartial}:
$$
\mu_{k,\ga'}(f*_\ga g)=\sum_{e=0}^k{{\mu_{e,\ga}(f)}
\over{[e]_{q}!}}{{q^{{(k-e)^2+k-e}\over2}[k]_{q}!}\over
{[k-e]_{q}!}}\int_{\ga'}g\,X^{k-e}=\sum_{e=0}^k\;
\qchoose k e \mu_{e,\ga}(f)\,\mu_{k-e,\ga'}(g).
$$
If $\ga=\ga'$ then,
by symmetry, the expression on the right
equals $\mu_{k,\ga}(g*_\ga f)=q^{{k^2+k}
\over2}\int_\ga (g*_{\ga}f)\,X^k$, hence the equality of the $q$-integrals.
The last statement is formula $(\ref{momentconvo})$ for $k=0$ and
$\ga=\ga'$.\qed
\begin{proposition}\label{convoleft}
Let $f\in \HD\II_{\ga,\al}^\om$
and $g\in \HD\II_{\ga',\be}^\om$.
Then $f*_{\ga}g\in \HD\II_{\ga',\eta}^\om$ with 
$\eta:={{\al\be}\over{\al+\be}}$.
\end{proposition}
\Proof By  Proposition \ref{fconvogintegrable},
$f*_{\ga}g\in\HD_a\II_{\ga'}$ for some $a>\ga'$
and 
by Proposition \ref{fede} $f*_{\ga}g\in\II_{\ga'}^\iy$.
It follows by equation (\ref{momentconvo}) that
\begin{equation}|\mu_{k,\ga'}(f*_\ga g)\bigr|\le  
\sum_{e=0}^k\,\qchoose k e
\bigl|\mu_{e,\ga}(f)\bigr|\,\bigl|\mu_{k-e,\ga'}(g)\bigr|.\end{equation}
By using the fact that
$\al e^2 +\beta (k-e)^2\ge k^2{{\al\be}\over{\al+\be}}$
we obtain that
\begin{equation}|\mu_{k,\ga'}(f*_\ga g)\bigr|\le C
q^{k^2{{\al\be}\over{2(\al+\be)}}}\sum_{e=0}^k b^ec^{k-e}\le C
q^{k^2{{\al\be}\over{2(\al+\be)}}}(\max(b,\,c))^k(k+1)\end{equation}
for some $b,c,C>0$.
\qed
\mLP
\begin{remark}
We have just seen that if $f\in \HD\II_{\ga,\al}^\om$
and $g\in \HD\II_{\ga',\be}^\om$ then
$f*_{\ga}g\in \HD\II_{\ga',\eta}^\om$,
where $\eta$ depends only on $\al$ and $\be$ and is
such that ${1\over{\eta}}={1\over\al}+{1\over\be}$.

\hfill$\spadesuit$
\end{remark}\smallbreak
We are ready to show that our convolution is associative.
\begin{theorem}\label{associativity}
Let $f\in \HD\II_{\ga}^\om$ and
$g\in \HD\II_{\ga'}^\om$. Let $h$ be defined,
together with its $q$-derivatives, on a domain $\Om$ and let 
$x\in\Om$ be such that, for some $R_x>0$,
$|(\pa^eh)(x)|=O(R_x^e)$ as $e\to\iy$. Then
$((f*_{\ga}g)*_{\ga'}h)(x)=(f*_{\ga}(g*_{\ga'}h))(x).$ 
In particular, the equality holds for every $x\in \Om$ if $h$ is
analytic on a starlike domain $\Om$ centered at $0$  and
it holds for every $x\in L(\ga'')$ if
$h\in\HD\II_{\ga''}^\iy$.  \end{theorem}
\Proof
By the previous results all series involved converge absolutely for $x$ as in
the hypothesis. We will show that the two given expressions coincide 
whenever they are well-defined. On the one hand
\begin{eqnarray*}
\lefteqn{((f*_{\ga}g)*_{\ga'}h)(x)=
\sum_{l=0}^\iy{{(-1)^l\mu_{l,\ga'}
(f*_\ga g)}\over{[\,l\,]_{q}!}}
(\pa^lh)(x)}\\
&&=\sum_{l=0}^\iy\sum_{k=0}^l{{(-1)^l
\mu_{k,\ga}(f)\mu_{l-k,\ga'}(g)}\over{[k]_q!\,[l-k]_q!}}(\pa^e h)(x)
\end{eqnarray*}
where we used  equation (\ref{momentconvo}) and Proposition \ref{fede}.
 
On the other hand
\begin{eqnarray*}
\lefteqn{(f*_\ga(g*_{\ga'}h))(x)=\sum_{l=0}^\iy{{(-1)^l
\mu_{l,\ga}(f)}\over{[l]_q!}}
\bigl(g*_{\ga'}\pa^lh\bigr)(x)}\hspace{.5in}\\
&&=\sum_{l=0}^\iy\sum_{m=0}^\iy{{(-1)^l\mu_{l,\ga}(f)}
\over{[l]_q!}}{{(-1)^m\mu_{m,\ga}(g)}\over{[m]_q!}}
(\pa^{l+m}h)(x).
\end{eqnarray*}
The two expressions will coincide if for one of them the double
sum is absolutely convergent.
This is certainly the case for $x$ satisfying $|\pa^eh(x)|=O(R_x^e)$
as $e\to\iy$. This last condition is in particular satisfied
if $h$ is analytic on a starlike domain $\Om$ centered at $0$ and
$x\in \Om$, 
or if $h\in\HD\II_{\ga''}^\iy$ and $x\in L(\ga'')$ (see also Lemma
\ref{interchange}).
\qed
\begin{corollary}\label{ideaal}
The class $\HD\II_\ga^\om$
is an algebra (not necessary unital) with respect to $*_{\ga}$.
Its subclass $\HS\II_\ga^\om$ is also an algebra (not necessary unital)
and it is a left ideal of $\HD\II_\ga^\om$.
\end{corollary}
\Proof By Theorem \ref{associativity} the convolution product in
$\HD\II_\ga^\om$ is   associative. $\HS\II_\ga^\om$ is a subclass of functions
of $\HD\II_\ga^\om$. By Lemma
\ref{welldefined} and Proposition \ref{convoleft} the product of a function in
$\HD\II_\ga^\om$ times a function in $\HS\II_\ga^\om$ is a function in
$\HD\II_\ga^\om$ which is analytic on a whole strip around $\bf R$, hence it
belongs to $\HS\II_\ga^\om$. By Theorem \ref{associativity} and the uniqueness
of analytic functions we have the statements.\qed
\smallbreak\smallbreak

We will see in Section 6, Example \ref{example1} that a left unit exists in
$\HS\II_\ga^\om$.
\section{Functions of strict left type}
Recall Definition \ref{left}(b) of the class $\II_\ga^\som$ of functions of
strict left type.
These functions can be characterized in the following way.
\begin{proposition}\label{schatting} Let $f\in\II_{\ga,\al}^\som$.
If $\al>1$ then 
$|f(\pm q^{-j}\ga)|=0$ for every $j\in\ZZ$;
if $\al=1$ then
$|f(\pm q^{-j}\ga)|=0$  for every $j\in\ZZ$ sufficiently large;
if $0<\al<1$ then, for some $c>0$,
$|f(\pm q^{-j}\ga)|=O(q^{{\be\over2} j^2}c^j)$ as $j\to\iy$
with $\be={1\over{1-\al}}$.

\noindent Conversely, let $f$ be a function on $L(\ga)$ which is bounded
on $\{\pm q^j\ga \mid j\ge1\}$ and let $\be>1$.
If, for some $c>0$,
$|f(\pm q^{-j}\ga)|=O(q^{{\be\over2} j^2}c^j)$ as $j\to\iy$
then $f\in\II_{\ga,\al}^\som$ with
$\al=1-\be^{-1}$. 
\end{proposition}
\Proof
Let $f\in\II_{\ga,\al}^\som$. Then there are
constants $C,b>0$ such that
$$
q^{{r^2+r}\over2}\sum_{k=-\iy}^\iy q^{k(1+r)}\ga^{r+1}\bigl[\,|f(
q^{k}\ga)|+|f(-q^{k}\ga)|\,\bigr]\le C q^{{\al\over2} r^2}b^r
\quad\hbox{for every  $r\in\Zplus$.}
$$
Hence
every term of the sum on the left-hand side is dominated by
the right-hand side. In particular,
\begin{equation}
|f(\pm q^{-j}\ga)|\le C
\,\ga^{-r-1}q^{{\al-1\over2}r^2+{2j-1\over2}r+j}b^r.
\label{bigger}\end{equation}
Therefore,
if $\al>1$, or if $\al=1$ and $q^{j-{1\over2}}<\ga b^{-1}$, we have
$|f(\pm q^{-j}\ga)|=0$ since the left hand side does not depend on $r$.
If
$0<\al<1$, since the left hand side of (\ref{bigger}) is independent of $r$, it
will be in particular smaller than the right-hand side evaluated for 
$r:=\bigl[{{2j-1}\over{2(1-\al)}}\bigr]\ge0$. Hence
\begin{eqnarray*}
\lefteqn{|f(\pm q^{-j}\ga)|\le C\ga^{-1}(b\ga^{-1})^r
q^{{\al-1\over2}r^2+{2j-1\over2}r+j}}\hspace{.5in}\\ 
&&\le C'
(b\ga^{-1})^{j\over 1-\al}q^{-{{(2j-1)^2(1-\al)}
\over{8(1-\al)^2}}
+{{(2j-1)^2}\over{4(1-\al)}}+{1\over2}}\\ 
&&\le
C'q^{1\over2}(b\ga^{-1})^{j\over 1-\al}q^{{(2j-1)^2}
\over{8(1-\al)}}=C''
\biggl((b\ga^{-1}q^{-{1\over2}})^{1\over1-\al}\biggr)^j
q^{{j^2}\over{2(1-\al)}},
\end{eqnarray*}
since
${{2j-1}\over{2(1-\al)}}-1\le r\le {{2j-1}\over{2(1-\al)}}$, and the
first statement follows.

\noindent
Next we will prove the converse statement. By boundedness of $f$ on the
set $\{\pm q^j\ga\mid j\ge0\}$ it follows that, for some $M>0$, we have
$\int_{-\ga}^{\ga}|t^kf(t)|d_{q}t\le 2M\ga^{k+1}$ for all $k\in\Zplus$. Then,
\begin{eqnarray*}
\lefteqn{q^{{k^2+k}\over2}\Biggl(\int_{-\gamma\cdot\infty}^{\ga\cdot\infty}
-\int_{-\gamma}^{\gamma}\Biggr)
|x^kf(x)| d_qx =(1-q)q^{{k^2+k}\over2}\sum_{j=1}^\iy
\sum_{\epsilon=\pm1}q^{-j(k+1)}\ga^{k+1}|f(\epsilon q^{-j}\ga)|}\\
&&\le 2\,C\,q^{{k^2+k}\over2}\ga^{k+1}\sum_{j=1}^\iy
q^{j^2{{\be}\over{2}}-j(k+1)}c^j\\
&&=2Cq^{-{1\over{2\be}}}\ga^{k+1}q^{{k\over2}-{k\over\be}}
q^{{k^2\over2}(1-\be^{-1})}\sum_{j=1}^\iy 
q^{{\be\over2}(j^2-2j{{(k+1)}\over\be}+{{(k+1)^2}
\over{\be^2}})}c^j\le M_2 b^k q^{{1\over2}\al k^2}
\end{eqnarray*}
for some  $M_2, b>0$ and with $\al:=1-\be^{-1}<1$, since the infinite sum is
dominated by
\sLP
$\displaystyle s(\max(c,1))^s +(\max(c,1))^{s+1}\sum_{j=0}^\iy
q^{{j^2}\over{2(1-\al)}}c^j,\;{\rm where}\;
s:=[(k+1)\be^{-1}]\ge0.$

\qed
\begin{remark}
By Proposition \ref{schatting} the pointwise product of functions
$f\in\II_{\ga,\al}^\som$ and $g\in\II_{\ga,\be}^\som$ satisfies the estimate
$|(fg)(\pm q^{-j}\ga)|\le C c^j q^{{{j^2}\over2}({1\over1-\al}+
{1\over1-\be})}$.
Hence $fg\in\II_{\ga,\eta}^\som$, where $\eta=1$ if
$\al$ or $\be=1$, and
$(1-\eta)^{-1}=(1-\al)^{-1}+(1-\be)^{-1}>2$ if $\al,\be<1$.
Thus $\eta>{1\over2}$ in all cases.
This will be useful for finding
examples in connection with commutativity of $q$-convolution (see Section
\ref{convocommu}).\hfill$\spadesuit$
\end{remark}
\begin{example}
As a consequence of Proposition \ref{schatting} all
functions belonging to $\HD\II_{\ga,1}^\som$
must be of the form $f\,E_{q^2}(-q^{2p}\ga^{-2}X^2)$ for some $p\in \ZZ$
and some analytic function $f$.\hfill$\spadesuit$ \end{example}
\begin{remark} Observe that functions of
strict left type $\al>1$ must be identically zero. From now on
$\HD\II^\som_{\ga,\al}$ will automatically imply
$\al\le1$.~{}\hfill$\spadesuit$
\end{remark}
\begin{example}
The family of functions $f_c$ ($c>0$) given by
$f_c(z):=e^{-c(\log(z^2+1))^2}$ provides for every $\al\in(0,\,1)$
an example of a function which is of strict left type $\al$ on $L(\ga)$
for each $\ga>0$.
Indeed, $f_c\in\HS_1$ and for every $\ga>0$ there exists $b>0$ such that
$|f_c(\pm q^{-k}\ga)|=O(q^{(4c\log(q^{-1}))k^2}b^k)$ as $k\to\iy$.
Hence $f\in\HS_1\II_{\ga,\al}^\som$ with
$\al:=1-(4c\log(q^{-1}))^{-1}\in(0,1)$ if
$c>\bigl((4\log(q^{-1}))^{-1}$.\hfill$\spadesuit$ \end{example}
\begin{example}
Let ${\cal E}_{q^2}(-X^2)$ be the function defined by
$$
{\cal E}_{q^2}(-x^2):=\sum_{k=0}^{\iy}{{(-1)^kq^{{k(k-1)}\over2}x^{2k}}
\over{(q^2;q^2)_k}}
={}_1\phi_1(0;-q;q,x^2).
$$
This function belongs with parameter value $1/2$
to a family of entire functions
interpolating between $e_{q^2}(-X^2)$ and $E_{q^2}(-X^2)$,
see also \cite{Ata}.
By formula (2.3) and Remark 2.4 in \cite{KooSw} with $z=-q$ and
$n=2k-1$ we have:
$$
(-q;q)_{\iy}\left|{\cal
E}_{q^2}(-q^{1-(2k-1)})\right|\le q^{2k^2-k}(-q;q)_{\iy}^2.
$$ 
Hence,
for $\ga=1$, $\left|{\cal E}_{q^2}(-(q^{-k}\ga)^2)\right|=O(q^{2k^2})$
as $k\to\iy$, which shows that
${\cal E}_{q^2}(-X^2)$ is of strict left type $3/4$ on
$L(1)$.  \hfill$\spadesuit$
\end{example}
\begin{lemma}\label{ciccia}
Let $f\in\HD\II_{\ga,\al}^\som$.
Then, there are $C,c,R>0$ such that 
$$
\nu_{e,\ga}(\pa^kf)\le C q^{{{\al}\over2}e^2}c^e R^k
$$
for all $k,\,e\in\Zplus$.
In particular, $\pa^kf\in\HD\II_{\ga,\al}^\som$.
\end{lemma}
\Proof
Let $f\in\HD_a\II_{\ga,\al}^\som$ and $r\in(\ga,\,a)$.
By Lemma \ref{integrapartial} there exists $B>0$ such that
$$
q^{{e^2+e}\over2}\int_\ga |X^{e}\,\pa^kf|\le
q^{{e^2+e}\over2}{{2^k}\over{\ga^k(1-q)^k}}
\Biggl[\int_\ga |X^{e}|\,|f|+
r^e B\Biggr]
$$
for all $k,\,e\in\Zplus$.
Combination with Definition \ref{left}(b) immediately yields the result.
\qed
\begin{proposition}\label{leftstrictstrict}
Let $f\in \HD\II_{\ga,\al}^\om$
and $g\in\HD\II_{\ga',\be}^\som$.
Then $f*_{\ga}g\in\HD\II_{\ga',\be}^\som$.
\end{proposition}
\Proof
For any $e\ge0$,
$$
q^{{e^2+e}\over2}\int_{\ga'}|f*_{\ga}g|\,|X^e|
\le
\sum_{k=0}^\iy{{|\mu_{k,\ga}(f)|}\over{[k]_q!}}q^{{e^2+e}\over2}
\int_{\ga'}|X|^e\,|\pa^kg|
\le C\sum_{k=0}^\iy{{q^{{\al\over2}k^2}b^k}\over{[k]_q!}}R^k
q^{{1\over2}\be e^2}c^e
$$
by formula (\ref{convoeq}), Definition \ref{left}(b)
and Lemma \ref{ciccia}.
\qed
\begin{corollary}\label{algebra}
The class $\HD\II_\ga^\som$  is a subalgebra 
of $\HD\II_\ga^\om$. Its subclass $\HS\II_\ga^\som$ is a
left ideal of $\HS\II_\ga^\om$, $\HD\II_\ga^\om$ and $\HD\II_\ga^\som$.
\end{corollary}
\Proof The proof is similar to that of Corollary \ref{ideaal}\hfill$\square$
\begin{definition}
For every $c\in[0,1)$ we denote by
$\HS\II_{\ga,>c}^\som$
the union of all spaces $\HS\II_{\ga,\al}^\som$ with $\al\in(c,1]$.\qed
\end{definition}
\begin{corollary}
The classes $\HS\II_{\ga,>c}^\som$ (for $c\in[0,1)$)
and $\HS\II_{\ga,c}^\som$ (for $c\in(0,1]$) are 
left ideals of $\HS\II_{\ga}^\som$ and of $\HS\II_{\ga}^\om$.
Similar properties hold for $\HD\II_{\ga}^\som$. \qed
\end{corollary}
\begin{corollary}
Let $f,g\in\HD\II_\ga^\om$, $h\in\HD$. Then
$(f*_{\ga}g*_{\ga}h)(x)=
(g*_{\ga}f*_{\ga}h)(x)$ for every $x$ where
the product is defined. In particular, 
for every pair of ideals $I\subset J$ with
$I,\,J\in\{\HS\II_\ga^\som,\,\HS\II_\ga^\om,\,
\HD\II_\ga^\om,\,\HD\II_\ga^\som\,\}$,
$I$ is a left module over $J/[J,\,J]_*$, where $[J,\,J]_*$
denotes the commutator ideal.
\end{corollary}
\Proof
This is a consequence of Lemma \ref{momenta} together with
the fact that $f*_{\ga}g$ depends only on the
$q$-moments $\mu_{e,\ga}(f)$  and not on the
values of $f(x)$.
The last statement follows by the inclusions
of left ideals
$\HS\II_\ga^\som\subset\HS\II_\ga^\om\subset\HD\II_\ga^\om$
and $\HS\II_\ga^\som\subset\HD\II_\ga^\som\subset\HD\II_\ga^\om$. \qed
\mPP
In particular, since $\HS\II_{\ga,>c}^\som\subset\HS\II_{\ga,>d}^\som$  
(respectively $\HS\II_{\ga,c}^\som\subset\HS\II_{\ga,d}^\som$) for
every $c\ge d$, we have a chain of left ideals on
$\HS\II_{\ga}^\som=\HS\II_{\ga,>0}^\som$, and similarly for 
$\HD\II_\ga^\som$. Hence
$E_{q^2}(-q^2X^2)\in\HS\II_{1,>c}^\som$ for every $c<1$ and
$e_{q^2}(-X^2)\in\HS\II_{\ga,>\half-\epsilon}^\som$ for every $\ga>0$
and $\epsilon\in(0,\,1/2)$.
\section{Commutativity of $q$-convolution}\label{convocommu}
We investigate commutativity now. We begin with a lemma.
\begin{lemma}\label{momentzero}
Let $f\in \HD\II_{\ga,>1/2}^\som$
be such that $\int_\ga f\,X^k=0$ for every $k\in\Zplus$. Then
$f(x)=0$ for every $x$ in some neighbourhood of zero.
In particular, if
$f\in \HS\II_{\ga,>1/2}^\som$, then $f(x)=0$ in each point $x$
where $f$ is analytic.
\end{lemma}
\Proof
Let $f$ have power series $f(x)=\sum_la_lx^l$ with
radius of convergence $>\ga$. Then
\begin{eqnarray*}
\lefteqn{\int_\ga |f|^2\,E_{q^2}(-q^2X^2\ga^{-2})
=\int_{-\ga}^\ga|f(x)|^2 E_{q^2}(-q^2x^2\ga^{-2})\,d_qx}\\
&&=\int_{-\ga}^\ga\biggl( f(x)E_{q^2}(-q^2x^2\ga^{-2})
\sum_{l=0}^\iy {\bar a}_lx^l\biggr)\,d_qx
=\sum_{l=0}^\iy {\bar a}_l\int_\ga f\,X^lE_{q^2}(-q^2\ga^{-2}X^2)\\
&&=\sum_{l=0}^\iy{\bar a}_l
\sum_{p=0}^\iy{{(-1)^p q^{p^2+p}\ga^{-2p}}\over{(q^2;q^2)_p}}
\int_\ga f\, X^{2p+l}
=0.
\end{eqnarray*}
Here the third and the fourth equality are justified by dominated
convergence. Indeed,
$$
\sum_{l=0}^{\iy}|a_l|
\int_{-\ga}^\ga|f( x)|\,|x^l|\, E_{q^2}(-q^2x^2\ga^{-2})\,d_{q}x\le
\sum_{l=0}^{\iy}|a_l|\ga^l\int_{-\ga}^{\ga}|f(x)|d_{q}x
<\iy
$$
and (use that $f$ has strict left type $>1/2$)
$$
\sum_{p=0}^\iy{{q^{p^2+p}\ga^{-2p}}\over{(q^2;q^2)_p}}\int_\ga 
|f\,X^{l+2p}|
\le 
C\sum_{p=0}^\iy{{q^{p^2+p}\ga^{-2p}}\over{(q^2;q^2)_p}}
q^{({1\over2}\al-{1\over2})(l+2p)^2}
(q^{-{1\over2}}b)^{l+2p}<\iy.
$$
Hence $|f(x)|^2=0$ if $x=\epsilon q^{k}\ga$ with $k\in\Zplus$ and
$\epsilon=\pm 1$.  
If moreover $f\in\HS\II_{\ga,>1/2}^\som$ then $f=0$
because it is analytic on a strip and vanishes on a sequence with
limit point in the strip.\qed
\mLP

\begin{remark}
Note the crucial role of analyticity of $f$ on a strip around
$\RR$ in order to conclude in the above Lemma that $f$ vanishes everywhere
on $L(\ga)$.

\hfill$\spadesuit$
\end{remark}
\begin{remark}
The proof of the above Lemma showed that for any
$f\in\HS$ 
the  power series of $f\,E_{q^2}(-q^{2}\ga^{-2}X^2)$ is absolutely
$q$-integrable on $L(\ga)$. \hfill$\spadesuit$
\end{remark}
\begin{theorem}\label{commutativity}
$\HS\II_{\ga,>c}^\som$ is a commutative algebra for every
$c\in[1/2,1)$. 
\end{theorem}
\Proof
If $f,g\in \HS\II_{\ga,>c}^\som$ then by Lemma \ref{momenta}, Proposition
\ref{leftstrictstrict} and Corollary \ref{algebra}
$f*_{\ga}g-g*_{\ga}f$ satisfies the hypothesis of Lemma \ref{momentzero}.\qed
\mPP
The following Theorem shows that $\HS_1\II_{\ga,1/2}^\som$ is
far from commutative as a $q$-convolution algebra.
Afterwards we give two other examples of noncommutativity.

\noindent Recall that $Q$ is the $q$-shift $Qf(x)=f(qx)$. 
\begin{theorem}\label{noncommut-ex}
Let $g(x):=e_{q^2}(-x^2)$ (so $g\in \HS_1\II_{\ga,1/2}^\som$ by Example
\ref{Gaussian} $\bf (a)$). Let $f\in\HS\II_{\ga,>1/2}^\som$ be an entire
function, not identically zero. Then
$(Q^{-n}f)*_\ga g\ne g*_\ga(Q^{-n}f)$
for $n$ sufficiently large.
\end{theorem}
\Proof
For the {\em discrete $q$-Hermite II polynomials}
$$
{\tilde h}_k(x;q)
:=
x^k{}_2\phi_1 (q^{-k},q^{-k+1};0;q^2,-q^2x^{-2})
=
(q;q)_k\sum_{l=0}^{[k/2]}{{(-1)^lq^{-2lk+2l^2+l}x^{k-2l}}
\over{(q^2;q^2)_l\,(q;q)_{k-2l}}}
$$
the following Rodrigues type formula was given in \cite{Koo},
formula (8.28): 
$$
(\pa^ke_{q^2}(-X^2))(x)={{(-1)^kq^{{k^2-k}\over2}}\over{(1-q)^k}}\,
{\tilde h}_k(x;q)\,e_{q^2}(-x^2).
$$
Also note by formula (II.6) in \cite{gasper} that
\begin{equation} \label{hki}
{\tilde h}_k(i;q)=i^k{}_2\phi_1 (q^{-k},q^{-k+1};0;q^2,q^2)
=
i^k q^{-\half k(k-1)}.
\end{equation}
For $f\in\II_\ga^\om$ put
\begin{equation}\label{f-tilde}
\tilde f(x):=\sum_{k=0}^\iy{q^{k^2-k\over 2}\mu_{k,\ga}(f)\over(q;q)_k}\,
\tilde h_k(x;q).
\end{equation}
It follows by the Rodrigues formula that
$$
(f*_{\ga}e_{q^2}(-X^2))(x)=\tilde f(x)\,e_{q^2}(-x^2)
$$
for $|\Im(x)|<1$.
In fact, $\tilde f$ is an entire function, which we will show by
uniform absolute convergence on compacta of the series defining
$\tilde f(x)$.
Indeed,
for $k=2h+\epsilon$ with $\epsilon=0,\,1$:
\begin{eqnarray*}
\lefteqn{{\tilde h}_k(x;q)=
(q;q)_{2h+\epsilon}\sum_{l=0}^{h}{{(-1)^lq^{-2l\epsilon
-4lh+2l^2+l}x^{2h+\epsilon-2l}}\over{(q^2;q^2)_l(q;q)_{2h+\epsilon-2l}}}}\\ 
&&=(-1)^h(q;q)_{2h+\epsilon}x^{\epsilon}q^{-2h^2-2h\epsilon+h}\sum_{p=0}^h
{{(-1)^pq^{2p^2-p+2p\epsilon}x^{2p}}
\over{(q^2;q^2)_{h-p}(q;q)_{2p}(1-q^{1+2p})^{\epsilon}}}.
\end{eqnarray*}
Hence
\begin{eqnarray*}\label{giovanna}
\lefteqn{{q^{k^2-k\over 2}\over(q;q)_k}\,|\tilde h_k(x;q)|
\le
{(q;q^2)_h(q^2;q^2)_h(1-q^{1+2h})^{\epsilon}|x|^{\epsilon}
\over(q;q)_k}}\hspace{1in}
\\
&&\times\sum_{p=0}^h
{{q^{2p^2-p}|x|^{2p}}
\over{(q^2;q^2)_{h}(q;q^2)_{h}(q^2;q^2)_p(1-q^{1+2p})^{\epsilon}}}
\le{\Phi(x)\over(q;q)_k}\,,\\ 
&&\qquad\qquad{\rm where}\quad
\Phi(x):=(1-q)^{-1}\max(1,|x|)\sum_{p=0}^\iy
{{q^{2p^2-p}|x|^{2p}}\over{(q^2;q^2)_p}}\,.
\end{eqnarray*}
In combination with Definition \ref{left}(a) this shows that,
for any $M>0$, the series
$\sum_{k=0}^\iy\Phi(x)\,{|\mu_{k,\ga}(f)|\over(q;q)_k}$ is uniformly
convergent in $x$ for $|x|\le M$.

Suppose that moreover $f\in\II_\ga^\som$.
Substitute the $q$-integral for $\mu_{k,\ga}(f)$ (see (\ref{momenteq}))
in equation (\ref{f-tilde}). Then we can
interchange $q$-integral and sum in the resulting expression,
by combination of the above estimates with  Definition \ref{left}(b).
It follows that $\tilde f$ can be seen as a $q$-integral transform
of $f$ with kernel $K(t,x)$: 
\begin{equation} \label{defkernel}
{\tilde f}(x)=
\int_{-\ga.\iy}^{\ga.\iy}f(t)\,K(t,x)d_qt,
\quad{\rm where}\quad
K(t,x):=\sum_{k=0}^\iy {{q^{k^2}t^k{\tilde h}_k(x;q)}\over{(q;q)_k}}\,.
\end{equation}
By equations (\ref{hki}) and (\ref{q-exponential}) we have
$K(t,i)=E_q(iqt)$.
Therefore 
$$
{\tilde f}(i)=\int_{-\ga.\iy}^{\ga.\iy}f(t)\,E_q(iqt)\,d_qt.
$$

{}From now on assume that $f\in\II_{\ga,>1/2}^\som$ and entire, not
identically zero.
For $\lambda\not=0$ put
$$
{\tilde f}_{\lambda}(x):=\lambda^{-1}
\int_{-\ga.\iy}^{\ga.\iy}f(\lambda^{-1}t)\,K(t,x)\,d_qt
=\int_{-\ga.\iy}^{\ga.\iy}f(t)\,K(\lambda t,x)\,d_qt.
$$
Then
\begin{equation}\label{transform}
{\tilde f}_{\lambda}(i)=
\int_{-\ga.\iy}^{\ga.\iy}f(t)\,E_q(iq\lambda t)\,d_qt=
\sum_{k=0}^\iy q^{-{k(k+1)\over 2}}\mu_{k,\ga}(f)\,(i\lambda)^k
\end{equation}
is well-defined for all $\lambda\in\CC$ and entire in $\lambda$.
If $\tilde f_\lambda(i)=0$ for $\lambda=q^{n_k}$, where
$n_k\to\iy$ in $\ZZ$ as $k\to\iy$, then $\tilde f_\lambda(i)=0$ for all
$\lambda\in\CC$.
Hence, by Lemma \ref{momentzero}, $f$ is identically zero, which
contradicts our assumption.
Thus for $n\in\ZZ$ sufficiently large we have
$\tilde f_{q^n}(i)\ne0$, hence
$(Q^{-n}f)*_\ga e_{q^2}(-X^2)=q^n\,e_{q^2}(-X^2)\,\tilde f_{q^n}$
does not extend to a function analytic at $i$. On the other hand
$e_{q^2}(-X^2)*_{\ga}(Q^{-n}f)$ is entire by
Lemma \ref{welldefined}. Hence the two products are different.
\qed
\begin{remark} \quad
\sLP
{\bf(a)}
Let $f\in\II_\ga^\som$. Then we can express the kernel $K(t,x)$ (defined by
(\ref{defkernel})) also as a $q$-hypergeometric function:
$$
K(t,x)
=(iqt;q)_\iy\,{}_1\phi_1(ix,iqt;q,-iqtx)
=(ix,-iqt;q)_\iy\,{}_2\phi_1(qtx^{-1},0;-iqt;q,ix).
$$
For the first identity use formula $(3.29.12)$ in \cite{KoSw}.
For the second identity use formula (III.1) in \cite{gasper}.
\mLP
{\bf (b)}
Since there are functions $f$ and $g$ in
$\HS\II_{\ga,1/2}^\som$
for which $F:=f*_{\ga}g-g*_{\ga}f\not\equiv0$, there exists a
function $F\in\HS\II_{\ga,1/2}^\som$, not identically zero,
for which $\int_\ga F\,X^k=0$ for every $k\in\Zplus$.
In view of Lemma \ref{momentzero} we can state that, for $\al\in(0,1]$,
the algebra $\HS\II_{\ga,\al}^\som$ is commutative iff each function $f$
in this algebra is determined by its moments $\mu_{e,\ga}(f)$ ($e\in\Zplus$).
\hfill$\spadesuit$
\end{remark}
\begin{example}\label{example1}
For $m\in\Zplus$ let
\begin{eqnarray}
\lefteqn{g_m(x):=
e_{q^2}(-x^2)\,{}_0\phi_1(-;q^{1+2m};q^2,-q^{1+2m}x^2)}
\hspace{.25in}\nonumber\\
&&=e_{q^2}(-x^2)
\sum_{r=0}^\iy{(-1)^r q^{2r(r-1)} q^{(1+2m)r} x^{2r} \over
(q^{1+2m};q^2)_r (q^2;q^2)_r}\,.\label{functions}
\end{eqnarray}
Then $g_m\in\HS_1$.
It will turn out that for each $\al>0$ the functions $g_0$ and $g_1$
are non-commuting elements of $\HS_1\II_{\ga,\al}^\om$.
First we consider more generally $g_m$.

By the fact that
$e_{q^2}(-X^2)\in\II_{\ga,1/2}^\som$ we can estimate that,
for some $b,C>0$,
\begin{eqnarray*}
\lefteqn{\int_\ga |g_m\,X^k|
\le\sum_{r=0}^\iy{q^{2r(r-1)}q^{(1+2m)r}\over(q^{1+2m};q^2)_r(q^2;q^2)_r}\,
q^{-{1\over2}(2r+k)(2r+k+1)}\nu_{2r+k,\ga}(e_{q^2}(-X^2))}
\hspace{1.5in}\\
&&\le C q^{-k^2-k\over2}b^k
\sum_{r=0}^\iy
{q^{{1\over4}(2r+k)^2}(bq^{m-1-k})^{2r}\over{(q^{1+2m};q^2)_r(q^2;q^2)_r}}<\iy.
\end{eqnarray*}
Hence $g_m\in\infigamma$ for every $\ga>0$.
By formula (\ref{eq-moment}), $\mu_{2k+1,\ga}(g_m)=0$ for every $k\in\Zplus$.
Furthermore, by dominated convergence we have
\begin{eqnarray*}
\lefteqn{\mu_{2k,\ga}(g_m)=
q^{2k^2+k}\sum_{r=0}^\iy{(-1)^r q^{2r(r-1)} q^{(1+2m)r}
\over{(q^{1+2m};q^2)_r(q^2;q^2)_r}}\int_\ga X^{2k+2r}
e_{q^2}(-X^2)}\hspace{.42in}\\ 
&&=c_{q}(\ga)q^{k^2+k}(q;q^2)_k
\sum_{r=0}^\iy{(q^{1+2k};q^2)_r (-1)^r q^{r^2-r} q^{(2m-2k)r}
\over{(q^{1+2m};q^2)_r (q^2;q^2)_r}}\\ 
&&=c_{q}(\ga)q^{k^2+k}(q;q^2)_k\,
{}_1\phi_1(q^{1+2k};q^{1+2m};q^2,q^{2m-2k})\\
&&=c_{q}(\ga)q^{k^2+k}(q;q^2)_k\,{{(q^{2m-2k};q^2)_{\iy}}
\over{(q^{1+2m};q^2)_{\iy}}}\,
\end{eqnarray*}
for every $k\in\Zplus$, where we used formula (\ref{eq-moment}) and 
formula (II.5) in \cite{gasper}.
Hence
$$
\mu_{e,\ga}(g_m)={c_q(\ga)\,(q^2;q^2)_\iy\over(q^{2m+1};q^2)_\iy}\,
{q^{k^2+k}(q;q^2)_k\over(q^2;q^2)_{m-k-1}}\quad
\hbox{if $e=2k$ with $k=0,1,\ldots,m-1$,}
$$
and $\mu_{e,\ga}(g_m)=0$ otherwise.
Hence $g_m\in\II_{\ga,\al}^\om$ for each $\al>0$ and
$$
(g_m*_\ga f)(x)={c_q(\ga)\,(q^2;q^2)_\iy\over(q^{2m+1};q^2)_\iy}\,
\sum_{k=0}^{m-1}{q^{k^2+k}\over(q^2;q^2)_k (q^2;q^2)_{m-k-1}}\,
(\pa^{2k}f)(x)
$$
for any function $f$ and for $x$ such that the $q$-derivatives of $f$ at
$x$ are well-defined.
In particular,
$$
g_0*_\ga f=0,\quad{\rm and}\quad
g_1*_\ga f={c_q(\ga)\,(q^2;q^2)_{\iy}
\over(q^3;q^2)_{\iy}}\,f\quad\hbox{for all $f$.}
$$
Hence,
$g_0*_{\ga}g_1=0\not=g_1*_{\ga}g_0$.

It follows from formulas (\ref{functions}) and (\ref{q-exponential}) that
$$
g_0(x)=e_{q^2}(-x^2)\,\Re(E_q(ix))=
\Re\left({(-ix;q)_\iy\over(-x^2;q^2)_\iy}\right)=
\Re\left({1\over(ix;q)_\iy}\right)=\Re(e_q(ix)).
$$
for $x\in\RR$. Hence
$$
g_0(\pm\ga q^{-k})=\Re\left({1\over(\pm i\ga q^{-k};q)_\iy}\right)=
\Re\left({(\pm i)^k\ga^{-k}q^{\half k(k+1)}\over
(\pm i\ga^{-1}q;q)_k(\pm i\ga;q)_\iy}\right)\quad
\hbox{for $k\in\Zplus$.}
$$
Therefore, the estimate
$|g_0(\pm\ga q^{-k})|=O(q^{{\be\over 2}k^2}c^k)$ as $k\to\iy$
is valid for some $c>0$ if and only if $\be\le 1$.
This result combined with Proposition \ref{schatting} implies that
$g_0$ is not of strict left type on $L(\ga)$.

The vanishing of all $q$-moments of $g_0$ can also be seen
directly from the formula
$$
\intiy x^n e_q(\pm ix)\,d_qx=
\intiy x^n e_{q^2}(-x^2)E_{q}(\pm ix)\,d_{q}x=0\quad
(n\in\Zplus,\;\ga>0),
$$
which follows from formula $(8.21)$ in
\cite{Koo} by substitution of $t=\pm q^{-1}$. Hence
$$
\qquad\qquad\qquad\intiy x^n g_0(x)\,d_qx=
\Re\intiy x^n e_q(ix)\,d_{q}x=0.\qquad\qquad\qquad\qquad\spadesuit
$$
\end{example}
\begin{example}\label{example2}
Define a function $g$ on $L(1)$ by
$$
g(\pm q^{k}):=(-1)^k q^{k} e_{q^2}(-q^{2k})\quad(k\in\ZZ).
$$
Then $g$
cannot be extended to a function in $\HD$ since it is alternating on a
sequence approaching to zero.
On the other hand,
$$
|g(\pm q^{-k})|={q^{-k}\over(-q^{-2k};q^2)_\iy}=
{q^{k^2}\over(-q^2;q^2)_k(-1;q^2)_\iy}=O(q^{k^2})\quad
\hbox{as $k\to\iy$,}
$$
and $|g(\pm q^k)|\le 1$ if $k\in\Zplus$.
Hence, by Proposition \ref{schatting} it follows that $g\in\II_{1,\half}^\som$.
Clearly, $\mu_{2n+1,1}(g)=0$ for all $n\in\Zplus$. Furthermore,
\begin{eqnarray*}
\lefteqn{\mu_{2n,1}(g)=2(1-q)\sum_{k=-\iy}^\iy (-1)^kq^{(2n+2)k}
e_{q^2}(-q^{2k})}\hspace{.5in}\\
&&=
{2(1-q)\over(-1;q^2)_\iy}\,{}_1\psi_1(-1,0;q^2,-q^{2n+2})=0
\end{eqnarray*}
for $n\in\Zplus$, where we used Ramanujan's ${}_1\psi_1$ summation
formula, see (II.29) in \cite{gasper}.
(C. Berg \cite{berg} used the same vanishing case of the ${}_1\psi_1$
in connection with the indeterminate moment problem related to
discrete $q$-Hermite II polynomials.)$\;$
We conclude that $g*_1f=0$ for every $f$.

Next we consider $f*_1g$ for some $f\in\II_1^\om$.
Since $g$ is not in $\HD$, we cannot use the results of Section
\ref{convoconve} in order to be sure that $(f*_1 g)(x)$ is
well-defined for suitable $x$.
However, we can reason as follows.
{}From the inequality
$|(\pa h)(x)|\le {|h(x)|+|h(qx)|\over(1-q)|x|}$ and the definition of $g$
we see by induction with respect to $e$ that
$$
|(\pa^eg)(\pm q^{k})|\le
e_{q^2}(-q^{2k}){{3^e q^{-k(e-1)}}\over{(1-q)^e}}
\quad(k,e\in\Zplus).
$$
Hence $(f*_1 g)(\pm q^k)$ is well-defined for $f\in\II_1^\om$,
$k\in\Zplus$.
Since $(-1)^kg(q^{k})>0$ for all $k\in\ZZ$,
one sees that $(-1)^k(\pa^eg)(q^{k})>0$ for all $k,e\in\Zplus$.
Hence $(f*_1g)(q^{2k})>0$ for $k\in\Zplus$ if $f\in\II_1^\om$ is
even and strictly positive (for instance $f:=e_{q^2}(-X^2)$).
Again we have obtained a couterexample to commutativity of
convolution.
\hfill$\spadesuit$
\end{example}
\mPP
Hence
we have shown by means of Theorem \ref{noncommut-ex}
and Examples \ref{example1}, \ref{example2}
that none of the hypotheses of Theorem \ref{commutativity}
($f,g\in\II_{\ga,>1/2}^\som$; $f\in\II_\ga^\som$; $f\in\HD$) can be relaxed.
\section{$q$-Convolution and $q$-Fourier transform}\label{convofourier}
In Remark \ref{fourier-remark} we introduced in (\ref{q-fourier})
a $q$-Fourier transform pair with the second transform
being given by
\begin{equation} \label{fouriereq1}
(\FF_\ga f)(y):=\int_{-\ga.\iy}^{\ga.\iy}E_q(iqxy)\,f(x)\,d_qx.
\end{equation}
The right-hand side of (\ref{fouriereq1}) can be formally rewritten
by power series expansion (\ref{q-exponential}) of $E_q(iqxy)$
and by substitution of formula (\ref{momenteq}) for
$\mu_{k,\ga}(f)$. This defines the following transform:
\begin{equation} \label{fouriereq2}
(\tilde\FF_\ga f)(y):=
\sum_{k=0}^\iy\mu_{k,\ga}(f)\,{(iy)^k\over(q;q)_k}\,.
\end{equation}
The transform $\tilde\FF_\ga$ is essentially the transform $F''_S(\id,\,\ga)$
in \cite{Ca2}  
for $n=1$, with the difference that in \cite{Ca2} $x$ and $y$ do not commute.
\begin{proposition}\quad\\
{\rm (a)}
If $f\in\II_\ga^\om$ then $\tilde\FF_\ga f$ is well-defined and it
is an entire analytic function.\\
{\rm (b)}
If moreover $f\in\II_\ga^\som$ then $\FF_\ga f$ is also well-defined
and $\FF_\ga f=\tilde\FF_\ga f$.\\
{\rm (c)}
Let $f\in\II_\ga^\om$. Then $\tilde\FF_\ga f=0$ iff $\mu_{k,\ga}(f)=0$ for all
$k\in\Zplus$.\\
{\rm (d)}
Let $f\in\HS\II_{\ga,>1/2}^\som$.
Then $\FF_\ga f=0$ iff $f=0$.\\
{\rm (e)}
Let $f\in\HD\II_\ga^\om$, $g\in\HD\II_{\ga'}^\om$.
Then $f*_\ga g\in\HD\II_{\ga'}^\om$ and
$\tilde \FF_\ga(f*_\ga g)=(\tilde \FF_\ga f)(\tilde\FF_{\ga'} g)$.
\end{proposition}
\Proof (a) follows by formula (\ref{fouriereq1}) and Definition \ref{left}(a).
\sLP
In order to prove (b), first substitute formula
(\ref{momenteq}) for $\mu_{k,\ga}(f)$ in formula (\ref{fouriereq2}).
Then justify interchange of summation and $q$-integration by dominated
convergence by using the
estimate for $\nu_{k,\ga}(f)$ in Definition \ref{left}(b).
Finally use the power series expansion (\ref{q-exponential}) for $E_q(iqxy)$. 
\sLP
(c) is evident from (a) together with (\ref{fouriereq1}).
\sLP
(d) follows from (c) and Lemma \ref{momentzero}.
\sLP
The first statement in (e) follows from Proposition \ref{convoleft}.
The second statement follows by taking
formula (\ref{fouriereq2}) for $f*_\ga g$ and then
substituting formula (\ref{momentconvo}).
Rearrangement of the double summation is justified by
dominated convergence (use Definition \ref{left}(a)).\qed
\section{Appendix}\label{appendix}
In our treatment of $f*_\ga g$ in this paper we usually required $g$
to be analytic on a neigbourhood of 0, while proofs of lemmas and
propositions only used that, for some $R>0$,
$|(\pa^k g)(x)|=O(R^k)$ as $k\to\iy$,
uniformly for $x\to 0$ in $L(\ga)$.
In this Appendix we show that the analyticity requirement on $g$ is not
an essential restriction.
In connection with the $k^{\rm th}$ $q$-derivative at 0 occurring
in the next Proposition, see also
Koekoek and Koekoek \cite{Koe}.
\begin{proposition}\label{analytic} Let $g$ be a function defined on the
half $q$-lattice $L^{\epsilon}(\ga)=\{\epsilon q^k\ga\}$ for
some $\ga>0$ and $\epsilon\in\{\pm1\}$. 
Suppose that
there exist constants $C>0$ and $r>\ga$ such that 
\begin{equation} \label{qder-estimate}
\bigl|(\pa^kg)(\epsilon q^t\ga)\bigr|\le
{{Cr}\over{(r-q^t\ga)}}\,{1\over{r^k(1-q)^k}}
\end{equation}
for every $k\ge0$ and for
every $t\in\ZZ$ for which $q^t\ga<r$.
Then the limit $l_p:=\lim_{k\to\iy}(\pa^pg)(\epsilon q^k\ga)$
exists and is finite for every $p\in\Zplus$, and
there exists a unique analytic function  $\tilde g$ on $\{x\in\CC\mid |x|<r\}$
such that
$\tilde g=g$ on $\{\epsilon q^k\ga\,|\, q^k\ga<r\}$.\\
If $g$ is defined  on the whole $q$-lattice $L(\ga)$ and if   
the condition above is satisfied for every $\epsilon\in\{\pm1\}$ and
if $l_p:=\lim_{k\to\iy}(\pa^pg)(\epsilon q^k\ga)$
is independent of $\epsilon$, then there exists a
unique analytic function  $\tilde g$ on $\{x\in\CC\mid |x|<r\}$ such that
$\tilde g=g$ on $\{\pm q^k\ga\mid q^k\ga<r\}$.\\
In particular, if $g$
is a function defined on $\bf R$ such that on every $q$-lattice $L(\ga)$ the
above conditions are satisfied and $l_p$ is independent of $\ga$
for every $p\in\Zplus$ 
then $g$ is analytic. 
\end{proposition}
\Proof
It follows from the identity
$$
g(\epsilon q^{k}\ga)=g(\epsilon q^{k+1}\ga)+
(1-q)\epsilon q^k\ga\,
(\pa g)(\epsilon q^{k}\ga)
$$
that
\begin{equation}\label{useful}
(\pa^bg)(\epsilon q^{k+p}\ga)=
(\pa^bg)(\epsilon
q^{k}\ga)-(1-q)\epsilon q^k\ga
\sum_{m=0}^{p-1}q^m\,(\pa^{b+1} g)(\epsilon q^{k+m}\ga)
\end{equation}
for every $b,p\in\Zplus$ and $k\in\ZZ$. 
For $k$ such that $q^k\ga<r$ it follows by (\ref{qder-estimate})
that
$\sum_{m=0}^\iy q^m\,(\pa^{b+1} g)(\epsilon q^{k+m}\ga)$
is absolutely convergent.
Hence both sides of equation (\ref{useful}) converge to a finite limit
$l_b$ as $p\to\iy$.

It follows by induction with respect to $n$ that
\begin{equation}
g(x)=\sum_{k=0}^n\,\qchoose n k (1-q)^k\,x^k\,
(\pa^k g)(q^{n-k}x).
\end{equation}
Let $x\in L(\ga)$, $|x|<r$. Then
\begin{eqnarray*}
\lefteqn{\qchoose n k (1-q)^k\,\left|x^k\,
(\pa^k g)(q^{n-k}x)\right|
\le {{(1-q)^k}\over{(q;q)_k}}\,|x|^k\,{{Cr}\over{r-q^{n-k}|x|}}\,
{1\over{r^k(1-q)^k}}}\hspace{2.5in}\\
&&\le
{{Cr}\over{r-|x|}}\,{{1}\over{(q;q)_k}}\,\biggl({{|x|}\over r}\biggr)^k.
\end{eqnarray*}
Hence by dominated convergence
we can take the limit for $n\to\iy$ and
$g(x)=\sum_{k=0}^\iy{{x^kl_k}\over{[k]_q!}}$
for $x=\epsilon q^k\ga$ with 
$q^k\ga<r$. \\
The other statements follow from the definition of
$\tilde g(x)=\sum_{k=0}^\iy{{x^kl_k}\over{[k]_q!}}$.
\qed
\begin{example}
Let $c_{q,\ga}(x)$ be defined on $L(\ga)$ as:
$$
c_{q,\ga}(\epsilon q^k\ga):=\cases{\int_{q^k\ga}e_{q^2}(-X^2)& if
$\epsilon=1$,\cr
 \int_{-q^k\ga}e_{q^2}(-X^2)& if $\epsilon=-1$.\cr}
$$
Then, since the $q$-integral is $q$-periodic,  the limit for $k\to\iy$ of
$c_{q,\ga}(\pm q^k\ga)=c_{q}(\ga)$, and
$(\pa^pc_{q,\ga})(x)\equiv 0$ on $L(\ga)$. The power series
$\tilde c_q(x)$ is then trivially the constant $c_q(\ga)$
defined by (\ref{ciqu}). However, the function
$c_q(z):=\int_{z}e_{q^2}(-X^2)$, coinciding with $c_{q,\ga}$ on every
$L(\ga)$ for $\ga>0$, is not analytic. Indeed the limit of $c_{q}(x)$ for
$x\to 0$ cannot exist. \hfill$\spadesuit$
\end{example}
\begin{example}
Consider the function $f(\epsilon q^k\ga):=(-\epsilon
q^k\ga;q)_{\iy}$ on  $L(\ga)$.
Clearly $\lim_{k\to\iy}f(\epsilon q^k\ga)=1$. One checks that 
$(\pa^nf)(\epsilon q^k\ga)={{q^{n(n-1)\over 2}}
\over{(1-q)^n}}(-\epsilon q^{k+n}\ga;q)_{\iy}$
so that
the limit for $k\to\iy$ is well-defined. The majorization for
$|(\pa^nf)(\epsilon q^k\ga)|$ is clearly verified for every 
$n\in\Zplus$, $k\in\ZZ$ and $\ga>0$.
Hence $f$ can be extended to the power series
$E_q(x)=\sum_{k=0}^\iy{{q^{{k(k-1)\over 2}}x^k}\over{(q;q)_k}}$. Since the
limits for $k\to\iy$ of the $q$-derivatives do not depend on $\ga$,
we have checked that 
$E_q(x)=(-x;q)_{\iy}$.

\hfill$\spadesuit$\end{example}
\vskip 0.7 truecm
\noindent
G. Carnovale,
Dipartimento di Matematica,\\
Universit\'a degli Studi di Roma,``Tor Vergata'',\\
via della Ricerca Scientifica 1,  00133 Roma, Italy;\\
email: {\tt carnoval@mat.uniroma2.it}
\bLP
T.~H. Koornwinder,
Korteweg-de Vries Instituut,\\
Universiteit van Amsterdam,\\
Plantage Muidergracht 24, 1018 TV Amsterdam,
The Netherlands;\\
email: {\tt thk@science.uva.nl}

\end{document}